\documentclass[a4paper,12pt]{article}
\usepackage{amsmath}
\usepackage{latexsym}
\usepackage{url}
\usepackage{color}
\usepackage{tikz}
\usetikzlibrary{intersections, calc, arrows.meta}
\numberwithin{equation}{section}

\def\R{{\bf R}}

\def\N{{\bf N}}

\def\d{\displaystyle}
\def\e{{\varepsilon}}

\def\p{\partial}
\def\v#1{\mbox{\boldmath $#1$}}

\newcommand{\LR}[1]{{\langle {#1} \rangle }}

\newtheorem{thm}{Theorem}[section]

\newtheorem{lem}{Lemma}[section]
\newtheorem{prop}{Proposition}[section]

\title{Semilinear wave equations of derivative type with characteristic weights\\
in one space dimension}
\author{
Shunsuke Kitamura
\footnote{
Doctor course, Mathematical Institute,
Tohoku University,
Aoba, Sendai 980-8578, Japan.
e-mail: shunsuke.kitamura.s8@dc.tohoku.ac.jp}
}
\date{
\[
\begin{array}{ll}
\mbox{\footnotesize{\bf Keywords:}}
& \mbox{\footnotesize semilinear wave equation, one dimension, classical solution, lifespan}\\
\mbox{\footnotesize{\bf MSC2020:}}
& \mbox{\footnotesize primary 35L71, secondary 35B44}\\
\end{array}
\]
}
\begin{document}
\maketitle
\begin{abstract}
In this paper, we investigate the lifespan estimates of classical solutions of the initial value problems for semilinear wave equations of derivative type with characteristic weights in one space dimension.
Such equations provide us basic principles on extending the general theory for nonlinear wave equations to the non-autonomous case.
In our results, two characteristic weights interact with each others
and produce a different range of parameters on the global-in-time existence
from the nonlinear terms of unknown function itself.
\end{abstract}


\section{Introduction}

\par 
We are concerned with the initial value problems;
\begin{equation}
\label{IVPderivative}
\left\{
\begin{array}{ll}
	\d u_{tt}-u_{xx}=\frac{|u_t|^p}{\LR{t+\LR{x}}^{1+a}\LR{t-\LR{x}}^{1+b}}
	&\mbox{in}\quad \R\times(0,\infty),\\
	u(x,0)=\e f(x),\ u_t(x,0)=\e g(x),
	& x\in\R,
\end{array}
\right.
\end{equation}
where $\LR{x} := \sqrt{1+x^2}$, $p>1$, $a,b \in\R$, $f$ and $g$ are given smooth functions of compact support
and a parameter $\e>0$ is \lq\lq small enough".
The future purpose is to extend the general theory for nonlinear wave equations
to the non-autonomous case.
Then, (\ref{IVPderivative}) may give us the principle in its setting on the nonlinear terms.
In this paper, we are interested in the estimate of the lifespan $T(\e)$,
the maximal existence time,
of classical solutions of (\ref{IVPderivative}). Our result is the following;
\begin{equation}\label{lifespan}
\begin{array}{l}
T(\e)=\infty\quad \mbox{for}\ a>0 \ \mbox{and} \ p(1+a)+b>0, \\
T(\e)\sim
\left\{
\begin{array}{lllll}
\exp(\e^{-(p-1)}) &\mbox{for}&  a=0\ \mbox{and}\ b\ge-p, \\
\exp(\e^{-p(p-1)}) &\mbox{for}&  a>0\ \mbox{and}\ p(1+a)+b=0, \\
\e^{-(p-1)/(-a)} &\mbox{for}& a<0 \ \mbox{and}\ b \ge -p,\\
\e^{-p(p-1)/(-p(1+a)-b)} &\mbox{for}& p(1+a)+b<0\ \mbox{and}\ b<-p.\\
\end{array}
\right.
\end{array}
\end{equation}
Here we defined $T(\e)\sim A(\e,C)$ as $A(\e,C_1)\le T(\e)\le A(\e,C_2)$ with positive constants,
$C_1$ and $C_2$, independent of $\e$.

When $a=b=-1$, the upper bounds in (\ref{lifespan}) are already obtained by Zhou \cite{Zhou01},
and the lower bounds are verified only for integer $p$ by general theory
which is studied by Li, Yu and Zhou \cite{LYZ91, LYZ92}. 
In addition, the initial value problem;
\[
\left\{
\begin{array}{ll}
	\d u_{tt}-u_{xx}=\frac{|u_t|^p}{\LR{x}^{1+a}}
	&\mbox{in}\quad \R\times(0,\infty),\\
	u(x,0)=\e f(x),\ u_t(x,0)=\e g(x),
	& x\in\R,
\end{array}
\right.
\]
which replaces the weights of (\ref{IVPderivative}) with the spatial weights, has the following results according to Kitamura, Morisawa and Takamura \cite{KMT22}:
\begin{equation}
\label{lifespan_x}
\begin{array}{l}
T(\e)=\infty\qquad\mbox{for}\quad a>0, \\
T(\e)\sim
\left\{
\begin{array}{ll}
\exp\left(C\e^{-(p-1)}\right) & \mbox{for}\ a=0, \\
C\e^{-(p-1)/(-a)} & \mbox{for}\ a<0,
\end{array}
\right.
\end{array}
\end{equation}
in which only the $C^1$ solution of the associated integral equation is considered for $1<p<2$.
But it can be also the classical solution by trivial modification on estimating
the nonlinear term with H\"older continuity.
(\ref{lifespan_x}) is consistent with (\ref{lifespan}) restricted by $b \ge -p$. 
The reason for this concidense is that lifespans are determined in the neighborhood of $t-|x|=0$ where all the quantities are equivalent, $\LR{t+\LR{x}} \sim \LR{x} \sim (1+t)$.

The case where the nonlinear term is a power of the unknown function itself with characteristic weights has already been considered in Kitamura, Takamura and Wakasa \cite{KTW22}. 
the lifespan estimates are classified into two cases according to
the value of the total integral of the initial speed $g$.
But (\ref{lifespan}) has no classification whatever it is.
This is due to the fact that Huygens' principle is always available
for the time derivative of the solution of the free wave equation.
In fact, 
\begin{equation}
\label{IVP}
\left\{
\begin{array}{ll}
	\d u_{tt}-u_{xx}=\frac{|u|^p}{\LR{t+\LR{x}}^{1+a}\LR{t-\LR{x}}^{1+b}}
	&\mbox{in}\quad \R\times(0,\infty),\\
	u(x,0)=\e f(x),\ u_t(x,0)=\e g(x),
	& x\in\R
\end{array}
\right.
\end{equation}
has the lifespan estimates: 
\begin{equation}
\label{lifespan_non-zero}
\begin{array}{ll}
T(\e)=\infty \quad \mbox{\rm if}\ a>0 \ \mbox{\rm and} \ a+b >0 \\
T(\e)\sim
\left\{
\begin{array}{lllll}
\exp(C\e^{-(p-1)}) &\mbox{\rm if}\ a+b=0\ \mbox{\rm and}\ a>0, \ \mbox{\rm or}\ a=0\ \mbox{\rm and} \ b>0, \\ 
\exp(C\e^{-(p-1)/2}) &\mbox{\rm if}\ a=b=0, \\
C\e^{-(p-1)/(-a)} &\mbox{\rm if}\ a<0 \ \mbox{\rm and}\ b>0,\\
\phi^{-1}(C\e^{-(p-1)})&\mbox{\rm if}\ a<0\ \mbox{\rm and}\ b=0,\\
C\e^{-(p-1)/(-a-b)}&\mbox{\rm if}\ a+b<0\ \mbox{\rm and}\ b<0,\\
\end{array}
\right.
\end{array}
\end{equation}
if 
\[
\int_\R g(x)dx\neq0,
\]
where $\phi^{-1}$ is an inverse function of $\phi$ defined by
\[
\phi(s):=s^{-a}\log(2+s),
\]
and
\begin{equation}
\label{lifespan_zero}
\begin{array}{ll}
T(\e)=\infty \quad \mbox{\rm if}\ a>0 \ \mbox{\rm and} \ a+b >0 \\
T(\e)\sim
\left\{
\begin{array}{lllll}
\exp(C\e^{-(p-1)}) &\mbox{\rm if}\  a=0\ \mbox{\rm and}\ b>0,\\
\exp(C\e^{-p(p-1)}) &\mbox{\rm if}\ a+b=0\ \mbox{\rm and} \ a>0, \\  
\exp(C\e^{-p(p-1)/(p+1)}) &\mbox{\rm if}\ a=b=0, \\
C\e^{-(p-1)/(-a)} &\mbox{\rm if}\ a<0 \ \mbox{\rm and}\ b>0,\\
\psi_1^{-1}(C\e^{-p(p-1)})&\mbox{\rm if}\ a<0\ \mbox{\rm and}\ b=0,\\
C\e^{-p(p-1)/(-pa-b)} &\mbox{\rm if}\ a<0\ \mbox{\rm and}\ b<0,\\
\psi_2^{-1}(C\e^{-p(p-1)})&\mbox{\rm if}\ a=0\ \mbox{\rm and}\ b<0,\\
C\e^{-p(p-1)/(-a-b)} &\mbox{\rm if}\ a+b<0 \ \mbox{\rm and}\ a>0,
\end{array}
\right.
\end{array}
\end{equation}
if 
\[
\int_\R g(x)dx=0,
\]
where $\psi_1^{-1},\psi_2^{-1}$ are inverse functions of $\psi_p$ defined by 
\[
\psi_1(s):=s^{-pa}\log(2+s), \ \psi_2(s):=s^{-b}\log^{p-1}(2+s).
\]
All the estimates of (\ref{lifespan}), (\ref{lifespan_non-zero}) and (\ref{lifespan_zero}) are represented on the $(a,b)$-plane as follows.

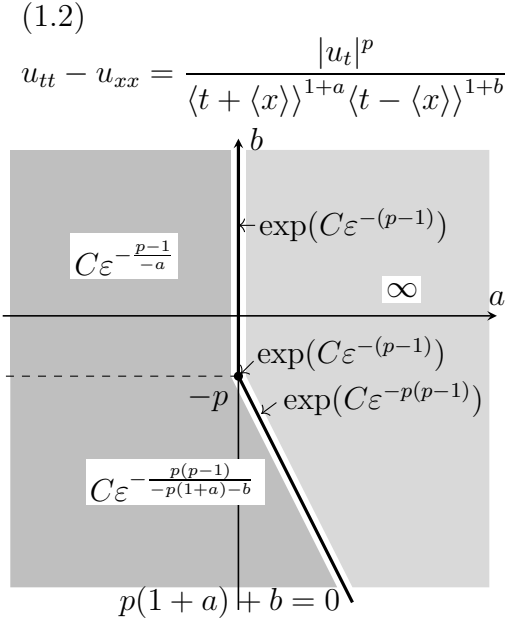
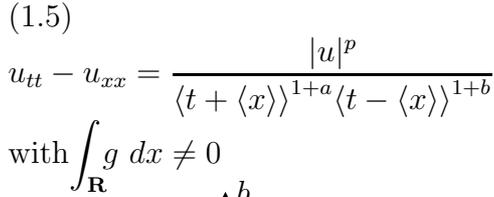
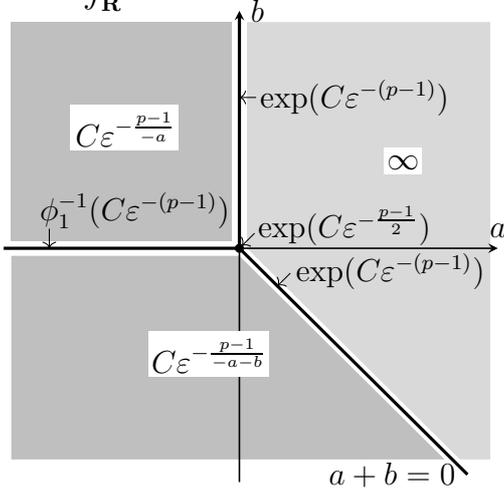
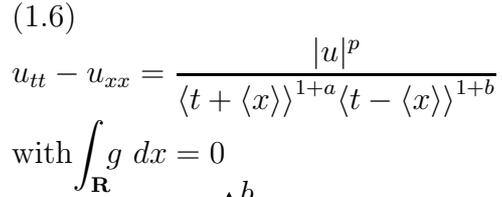
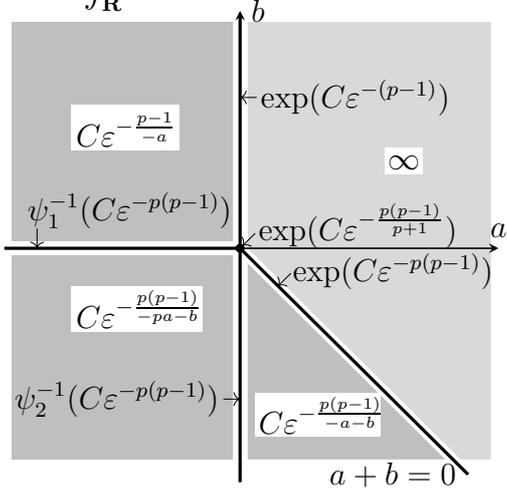
\centerline{
\begin{tikzpicture} 
 \fill[white!40!lightgray](0.1,0)--(1.5,-2.8)--(3.3,-2.8)--(3.3,3)--(0.1,3)--cycle;
 \fill[lightgray](-0.1,0)--(-3,0)--(-3,3)--(-0.1,3)--cycle;
 \fill[lightgray](-0.1,0)--(-3,0)--(-3,-2.8)--(1.3,-2.8)--cycle;
 \draw[->,>=stealth,semithick] (-3.1,0.8)--(3.4,0.8)node[above]{$a$}; 
 \draw[->,>=stealth,semithick] (0,-3.1)--(0,3.15)node[right]{$b$}; 
 \draw (-3,4.2)node[right,align=left]{(\ref{lifespan})\\$\d u_{tt}-u_{xx} = \frac{|u_t|^p}{\LR{t+\LR{x}}^{1+a}\LR{t-\LR{x}}^{1+b}}$};
 \fill [white] (1.9,0.98) rectangle (2.4,1.33);
 \draw (1.8,0.9)node[above right]{$\infty$}; 
 \draw[->] (0.5,-0.35)--(0.3,-0.5);
 \draw (0.45,-0.3)node[right]{\rm exp$(C\e^{-p(p-1)})$};
 \draw[->] (0.2,0.2)--(0.03,0.03);
 \draw (0.1,0.3)node[right]{\rm exp$(C\e^{-(p-1)})$};
 \draw[->] (0.22,2)--(0,2);
 \draw (0.13,2)node[right]{\rm exp$(C\e^{-(p-1)})$};
 \fill [white] (-0.8,1.3) rectangle (-2.22,1.9);
 \draw (-2.3,1.2)node[above right]{$C\e^{-\frac{p-1}{-a}}$}; 
 \draw[dashed] (0,0)--(-3.1,0);
 \draw (-0.4,-0.3)node{$-p$};
 \fill [white] (-2,-1.7) rectangle (0.35,-1.1);
 \draw (-2.1,-1.8)node[above right]{$C\e^{-\frac{p(p-1)}{-p(1+a)-b}}$};
 \draw[very thick,domain=0:3.1] plot(0,\x);
 \draw[very thick,domain=0:1.5] plot(\x,-2*\x)node[left]{$p(1+a)+b=0$};
  \fill [black] (0,0) circle [radius=0.06];
\end{tikzpicture}
}

\noindent
\begin{tikzpicture} 
 \fill[white!40!lightgray](0.1,0)--(2.9,-2.8)--(3.3,-2.8)--(3.3,3)--(0.1,3)--cycle;
 \fill[lightgray](-0.1,0.1)--(-3,0.1)--(-3,3)--(-0.1,3)--cycle;
 \fill[lightgray](0,-0.1)--(-3,-0.1)--(-3,-2.8)--(2.7,-2.8)--cycle;
 \draw[->,>=stealth,semithick] (-3.1,0)--(3.4,0)node[above]{$a$};
 \draw[->,>=stealth,semithick] (0,-3.1)--(0,3.15)node[right]{$b$};
 \draw (-3,4.4)node[right,align=left]{(\ref{lifespan_non-zero})\\$\d u_{tt}-u_{xx} = \frac{|u|^p}{\LR{t+\LR{x}}^{1+a}\LR{t-\LR{x}}^{1+b}}$ \\ with$\d \int_\R \!\! g \ dx \neq 0$};
 \fill [white] (1.9,0.98) rectangle (2.4,1.33);
 \draw (1.8,0.9)node[above right]{$\infty$}; 
 \draw[->] (0.67,-0.33)--(0.5,-0.5);
 \draw (0.6,-0.3)node[right]{\rm exp$(C\e^{-(p-1)})$};
 \draw[->] (0.2,0.2)--(0.03,0.03);
 \draw (0.1,0.3)node[right]{\rm exp$(C\e^{-\frac{p-1}{2}})$};
 \draw[->] (0.22,2)--(0,2);
 \draw (0.13,2)node[right]{\rm exp$(C\e^{-(p-1)})$};
 \fill [white] (-0.8,1.3) rectangle (-2.22,1.9);
 \draw (-2.3,1.2)node[above right]{$C\e^{-\frac{p-1}{-a}}$}; 
 \draw[->] (-2.5,0.26)--(-2.5,0);
 \draw (-2.77,0.5)node[right]{$\phi_{1}^{-1}(C\e^{-(p-1)})$};
 \fill [white] (-1.2,-1.7) rectangle (0.38,-1.1);
 \draw (-1.3,-1.8)node[above right]{$C\e^{-\frac{p-1}{-a-b}}$};
 \draw[very thick,domain=-3.1:0] plot(\x,0);
 \draw[very thick,domain=0:3.1] plot(0,\x);
 \draw[very thick,domain=0:3] plot(\x,-\x)node[left]{$a+b=0$};
   \fill [black] (0,0) circle [radius=0.06];
\end{tikzpicture}
\begin{tikzpicture} 
 \fill[white!40!lightgray](0.1,0)--(2.9,-2.8)--(3.3,-2.8)--(3.3,3)--(0.1,3)--cycle;
 \fill[lightgray](-0.1,0.1)--(-3,0.1)--(-3,3)--(-0.1,3)--cycle;
 \fill[lightgray](-0.1,-0.1)--(-3,-0.1)--(-3,-2.8)--(-0.1,-2.8)--cycle;
 \fill[lightgray](0.1,-0.2)--(0.1,-2.8)--(2.7,-2.8)--cycle;
 \draw[->,>=stealth,semithick] (-3.1,0)--(3.4,0)node[above]{$a$};
 \draw[->,>=stealth,semithick] (0,-3.1)--(0,3.15)node[right]{$b$};
 \draw (-3,4.4)node[right,align=left]{(\ref{lifespan_zero})\\$\d u_{tt}-u_{xx} = \frac{|u|^p}{\LR{t+\LR{x}}^{1+a}\LR{t-\LR{x}}^{1+b}}$ \\ with$\d \int_\R \!\! g \ dx =0$};
 \fill [white] (1.9,0.98) rectangle (2.4,1.33);
 \draw (1.8,0.9)node[above right]{$\infty$}; 
 \draw[->] (0.67,-0.33)--(0.5,-0.5);
 \draw (0.55,-0.3)node[right]{\rm exp$(C\e^{-p(p-1)})$};
 \draw[->] (0.2,0.2)--(0.03,0.03);
 \draw (0.1,0.3)node[right]{\rm exp$(C\e^{-\frac{p(p-1)}{p+1}})$};
 \draw[->] (0.22,2)--(0,2);
 \draw (0.13,2)node[right]{\rm exp$(C\e^{-(p-1)})$};
 \fill [white] (-0.8,1.3) rectangle (-2.22,1.9);
 \draw (-2.3,1.2)node[above right]{$C\e^{-\frac{p-1}{-a}}$}; 
 \draw[->] (-2.67,0.26)--(-2.67,0);
 \draw (-2.94,0.5)node[right]{$\psi_1^{-1}(C\e^{-p(p-1)})$};
 \fill [white] (-0.5,-0.5) rectangle (-2.22,-1.1);
 \draw (-2.3,-1.2)node[above right]{$C\e^{-\frac{p(p-1)}{-pa-b}}$}; 
 \draw[->] (-0.22,-2)--(0,-2);
 \draw (-0.13,-2)node[left]{$\psi_2^{-1}(C\e^{-p(p-1)})$};
 \fill [white] (0.2,-2.5) rectangle (1.85,-1.92);
 \draw (0.1,-2.6)node[above right]{$C\e^{-\frac{p(p-1)}{-a-b}}$}; 
 \draw[very thick,domain=-3.1:0] plot(\x,0);
 \draw[very thick,domain=-3.1:3.1] plot(0,\x);
 \draw[very thick,domain=0:3] plot(\x,-\x)node[left]{$a+b=0$};
   \fill [black] (0,0) circle [radius=0.06];
\end{tikzpicture}
We remark that the quantities in all the cases of (\ref{lifespan_non-zero}) are smaller than
those of (\ref{lifespan_zero}), and those of (\ref{lifespan_zero}) are smaller than those of (\ref{lifespan}). In the fourth quadrant where the characteristic weights interact, the lines of the case of almost global which means lifespan grows exponentially, are different for $|u_t|^p$ and $|u|^p$. The reason why the slope of the line of (\ref{lifespan}) is $p$ is that its optimal lifespan is obtained by the second iteration, not the first iteration. 


\par
This paper is organized as follows.
In the next section, (\ref{lifespan}) is divided into two theorems,
and the preliminaries are introduced.
Section 3 is devoted to the proof of the existence part of (\ref{lifespan}).
The main strategy is the iteration method in the weighted $L^\infty$ space 
which is originally introduced by John \cite{John79}.
In Section 4, we prove a priori estiamte.
Finally, we prove the blow-up part of (\ref{lifespan}) by two propositions. 
The first proposition is proved employing the method by Zhou \cite{Zhou01} in Section 5, 
and the second proposition is proved by the iteration method in the weighted $L^\infty$ space, similar to the method of Kitamura, Takamura and Wakasa \cite{KTW22} in Section 5.


\section{Preliminaries and main results}

Throughout of this paper, we assume that the initial data
$(f,g)\in C_0^2(\R)\times C^1_0(\R)$ satisfies
\begin{equation}
\label{supp_initial}
\mbox{\rm supp }f,\ \mbox{supp }g\subset\{x\in\R:|x|\le R\},\quad R\ge1.
\end{equation}
Let $u$ be a classical solution of (\ref{IVPderivative}) in the time interval $[0,T]$.
Then the support condition of the initial data, (\ref{supp_initial}), implies that
\begin{equation}
\label{support_sol}
\mbox{supp}\ u(x,t)\subset\{(x,t)\in\R\times[0,T]:|x|\le t+R\}.
\end{equation}
For example, see Appendix in John \cite{John_book} for this fact.

\par
It is well-known that $u$ satisfies the following integral equation;
\begin{equation}
\label{integral}
u(x,t)=\e u^0(x,t)+L_{a,b}(|u_t|^p)(x,t),
\end{equation}
where $u^0$ is a solution of the free wave equation with the same initial data;
\begin{equation}
\label{linear}
u^0(x,t):=\frac{1}{2}\{f(x+t)+f(x-t)\}+\frac{1}{2}\int_{x-t}^{x+t}g(y)dy,
\end{equation}
and a linear integral operator $L_{a,b}$ for a function $v=v(x,t)$ is Duhamel's term defined by
\begin{equation}
\label{nonlinear}
L_{a,b}(v)(x,t):=\frac{1}{2}\int_0^tds\int_{x-t+s}^{x+t-s}\frac{v(y,s)}{\LR{s+\LR{y}}^{1+a}\LR{s-\LR{y}}^{1+b}}dy.
\end{equation}
Then, one can apply the time-derivative to (\ref{integral}) and (\ref{linear}) to obtain
\begin{equation}
\label{integral_derivative}
u_t(x,t)=\e u_t^0(x,t)+L_{a,b}'(|u_t|^p)(x,t)
\end{equation}
and
\begin{equation}
\label{linear_derivative}
u_t^0(x,t)=\frac{1}{2}\{f'(x+t)-f'(x-t)+g(x+t)+g(x-t)\},
\end{equation}
where $L_{a,b}'$ for a function $v=v(x,t)$ is defined by
\begin{equation}
\label{nonlinear_derivative}
\begin{array}{ll}
L_{a,b}'(v)(x,t)
&\d:=\frac{1}{2}\int_0^t\frac{v(x+t-s,s)}{\LR{s+\LR{x+t-s}}^{1+a}\LR{s-\LR{x+t-s}}^{1+b}}ds\\
&\quad+\d\frac{1}{2}\int_0^t\frac{v(x-t+s,s)}{\LR{s+\LR{x-t+s}}^{1+a}\LR{s-\LR{x-t+s}}^{1+b}}ds.
\end{array}
\end{equation}
On the other hand, applying the space-derivative to (\ref{integral}) and (\ref{linear}),
we have
\[
u_x(x,t)=\e u_x^0(x,t)+\overline{L_a'}(|u_t|^p)(x,t)
\]
and
\[
u_x^0(x,t)=\frac{1}{2}\{f'(x+t)+f'(x-t)+g(x+t)-g(x-t)\},
\]
where $\overline{L_{a,b}'}$ for a function $v=v(x,t)$ is defined by
\begin{equation}
\label{nonlinear_derivative_conjugate}
\begin{array}{ll}
\overline{L_{a,b}'}(v)(x,t):=
&\d\frac{1}{2}\int_0^t\frac{v(x+t-s,s)}{\LR{s+\LR{x+t-s}}^{1+a}\LR{s-\LR{x+t-s}}^{1+b}}ds\\
&-\d\frac{1}{2}\int_0^t\frac{v(x-t+s,s)}{\LR{s+\LR{x-t+s}}^{1+a}\LR{s-\LR{x-t+s}}^{1+b}}ds.
\end{array}
\end{equation}
Therefore, $u_x$ is expressed by $u_t$.
Moreover, one more space-derivative to (\ref{integral_derivative}) yields that
\begin{equation}
\label{integral_2derivative}
\begin{array}{ll}
u_{tx}(x,t)=\e u_{tx}^0(x,t)&\d+pL_{a,b}'(|u_t|^{p-2}u_tu_{tx})(x,t) \\
&\d -(1+a)L_{a+1,b}'\left(|u_t|^p\frac{x}{\LR{x}}\frac{t+\LR{x}}{\LR{t+\LR{x}}}\right)(x,t)\\
&\d +(1+b)L_{a,b+1}'\left(|u_t|^p\frac{x}{\LR{x}}\frac{t-\LR{x}}{\LR{t-\LR{x}}}\right)(x,t)
\end{array}
\end{equation}
and
\begin{equation}
\label{linear_2derivative}
u_{tx}^0(x,t):=\frac{1}{2}\{f''(x+t)-f''(x-t)+g'(x+t)+g'(x-t)\}.
\end{equation}
because 
\[
\begin{array}{ll}
\d \frac{d \LR{t+\LR{x}}^{-1-a}}{dx} 
&\d= \frac{d (t+\LR{x})}{dx} \frac{d \LR{t+\LR{x}}}{d(t+\LR{x})} \frac{d \LR{t+\LR{x}}^{-1-a}}{d\LR{t+\LR{x}}} \\
&\d= -(1+a) \frac{x}{\LR{x}} \frac{t+\LR{x}}{\LR{t+\LR{x}}} \frac{1}{\LR{t+\LR{x}}^{2+a}}
\end{array}
\]
holds.
Similarly, we have that
\begin{equation}
\label{u_{tt}}
\begin{array}{ll}
u_{tt}(x,t)=&
\d\e u_{tt}^0(x,t)+\frac{|u_t(x,t)|^p}{\LR{t+\LR{x}}^{1+a}\LR{t-\LR{x}}^{1+b}} +p\overline{L_{a,b}'}(|u_t|^{p-2}u_tu_{tx})(x,t) \\
&\d -(1+a)\overline{L_{a+1,b}'}\left(|u_t|^p\frac{x}{\LR{x}}\frac{t+\LR{x}}{\LR{t+\LR{x}}}\right)(x,t)\\
&\d +(1+b)\overline{L_{a,b+1}'}\left(|u_t|^p\frac{x}{\LR{x}}\frac{t-\LR{x}}{\LR{t-\LR{x}}}\right)(x,t)
\end{array}
\end{equation}
and
\[
u_{tt}^0(x,t)=\frac{1}{2}\{f''(x+t)+f''(x-t)+g'(x+t)-g'(x-t)\}.
\]
Therefore, $u_{tt}$ is expressed by $u_{tx}$ and $u_t$, so is $u_{xx}$ because of
\[
\begin{array}{ll}
u_{xx}(x,t)=
&\d\e u_{xx}^0(x,t)+p\overline{L_{a,b}'}(|u_t|^{p-2}u_tu_{tx})(x,t) \\
&\d -(1+a)\overline{L_{a+1,b}'}\left(|u_t|^p\frac{x}{\LR{x}}\frac{t+\LR{x}}{\LR{t+\LR{x}}}\right)(x,t)\\
&\d +(1+b)\overline{L_{a,b+1}'}\left(|u_t|^p\frac{x}{\LR{x}}\frac{t-\LR{x}}{\LR{t-\LR{x}}}\right)(x,t)
\end{array}
\]
and
\[
u_{xx}^0(x,t)=u^0_{tt}(x,t).
\]

\par
First, we note the following fact.

\begin{prop}
\label{prop:continuity}
Assume that $(f,g)\in C^2(\R)\times C^1(\R)$.
Let $u_t$ be a $C^1$ solution of (\ref{integral_derivative}). 
Then,
\begin{equation}
\label{set}
w(x,t):=\int_0^tu_t(x,s)ds+\e f(x)
\end{equation}
 is a classical solution of (\ref{IVPderivative}).
\end{prop}
\par\noindent
{\bf Proof.} It is trivial that $w$ satisfies the initial condition and
\begin{equation}
\label{equality}
w_t=u_t,\quad w_{tt}=u_{tt}.
\end{equation}
Then, (\ref{integral_2derivative}) yields that
\[
\begin{array}{ll}
&\d w_x(x,t) \\
&\d=\int_0^tu_{tx}(x,s)ds+\e f'(x)\\
&\d=\int_0^tpL_{a,b}'(|u_t|^{p-2}u_tu_{tx})(x,s)-(1+a)L_{a+1,b}'\left(|u_t|^p\frac{x}{\LR{x}}\frac{t+\LR{x}}{\LR{t+\LR{x}}}\right)(x,s)\\
&\d\quad+(1+b)L_{a,b+1}'\left(|u_t|^p\frac{x}{\LR{x}}\frac{t-\LR{x}}{\LR{t-\LR{x}}}\right)(x,s)ds+\int_0^t\e u_{tx}^0(x,s)ds+\e f'(x)\\
&\d=\overline{L'_{a,b}}(|u_t|^p)(x,t)+\e u_x^0(x,t)
\end{array}
\]
because of
\[
\begin{array}{ll}
\d pL_{a,b}'(|u_t|^{p-2}u_tu_{tx})(x,s)\\
\d -(1+a)L_{a+1,b}'\left(|u_t|^p\frac{x}{\LR{x}}\frac{t+\LR{x}}{\LR{t+\LR{x}}}\right)(x,s)\\
\d +(1+b)L_{a,b+1}'\left(|u_t|^p\frac{x}{\LR{x}}\frac{t-\LR{x}}{\LR{t-\LR{x}}}\right)(x,s)\\
\d =\frac{\p}{\p s}\overline{L'_{a,b}}(|u_t|^p)(x,s).
\end{array}
\]
Therefore we obtain that
\[
\begin{array}{ll}
w_{xx}(x,t)=
&\d\e u_{xx}^0(x,t)+p\overline{L_{a,b}'}(|u_t|^{p-2}u_tu_{tx})(x,t) \\
&\d -(1+a)\overline{L_{a+1,b}'}\left(|u_t|^p\frac{x}{\LR{x}}\frac{t+\LR{x}}{\LR{t+\LR{x}}}\right)(x,t)\\
&\d +(1+b)\overline{L_{a,b+1}'}\left(|u_t|^p\frac{x}{\LR{x}}\frac{t-\LR{x}}{\LR{t-\LR{x}}}\right)(x,t)
\end{array}
\]
which implies, together with (\ref{u_{tt}}) and (\ref{equality}), the desired conclusion,
\[
w_{tt}-w_{xx}=\frac{|u_t|^p}{\LR{t+\LR{x}}^{1+a}\LR{t-\LR{x}}^{1+b}}=\frac{|w_t|^p}{\LR{t+\LR{x}}^{1+a}\LR{t-\LR{x}}^{1+b}}.
\]
\hfill$\Box$

\vskip10pt
Our result in (\ref{lifespan}) is splitted into the following two theorems.
\begin{thm}
\label{thm:lower-bound}
Assume (\ref{supp_initial}).
Then, there exists a positive constant $\e_1=\e_1(f,g,p,a,b,R)>0$ such that
a classical solution $u\in C^2(\R\times[0,T])$ of (\ref{IVPderivative})
exists as far as $T$ satisfies
\begin{equation}
\label{lower-bound}
\begin{array}{l}
T\le
\left\{
\begin{array}{lllll}
\exp(c\e^{-(p-1)}) &\mbox{for}&  a=0\ \mbox{and}\ b\ge-p, \\
\exp(c\e^{-p(p-1)}) &\mbox{for}&  a>0\ \mbox{and}\ p(1+a)+b=0, \\
c\e^{-(p-1)/(-a)} &\mbox{for}& a<0 \ \mbox{and}\ b \ge -p,\\
c\e^{-p(p-1)/(-p(1+a)-b)} &\mbox{for}& p(1+a)+b<0\ \mbox{and}\ b<-p.\\
\end{array}
\right.
\\
T<\infty\quad \mbox{for}\ a>0 \ \mbox{and} \ p(1+a)+b>0,
\end{array}
\end{equation}
where $0<\e\le\e_1$, $c$ is a positive constant independent of $\e$.
\end{thm}
\begin{thm}
\label{thm:upper-bound}
Assume (\ref{supp_initial})
and
\begin{equation}
\label{positive_non-zero}
\int_{\R}g(x)dx>0.
\end{equation}
Then, there exists a positive constant $\e_2=\e_2(g,p,a,b,R)>0$ such that
a solution $u_t\in C(\R\times[0,T])$ with
\[
\mbox{\rm supp}\ u_t(x,t)\subset\{(x,t)\in\R\times[0,T]:|x|\le t+R\}
\]
of associated integral equations (\ref{integral_derivative})
to (\ref{IVPderivative})
cannot exist whenever $T$ satisfies
\begin{equation}
\label{upper-bound}
T\ge
\left\{
\begin{array}{lllll}
\exp(C\e^{-(p-1)}) &\mbox{for}&  a=0\ \mbox{and}\ b\ge-p, \\
\exp(C\e^{-p(p-1)}) &\mbox{for}&  a>0\ \mbox{and}\ p(1+a)+b=0, \\
C\e^{-(p-1)/(-a)} &\mbox{for}& a<0 \ \mbox{and}\ b \ge -p,\\
C\e^{-p(p-1)/(-p(1+a)-b)} &\mbox{for}& p(1+a)+b<0\ \mbox{and}\ b<-p.\\
\end{array}
\right.
\end{equation}
where $0<\e\le\e_2$, $C$ is a positive constant independent of $\e$.
\end{thm}
The proofs of above theorems are given in following sections.


\section{Proof of Theorem \ref{thm:lower-bound}}
\par
According to the observations in the previous section,
we shall construct a $C^1$ solution of (\ref{integral_derivative}). 
Let $\{U_j(x,t)\}_{j\in\N}$ be a sequence of $C^1(\R\times[0,T])$ defined by
\begin{equation}
\label{U_j}
U_{j+1}=L'_{a,b}(|U_j+\e u_t^0|^p),\ U_1=0.
\end{equation}
Then, in view of (\ref{integral_2derivative}), $(U_j)_x$ has to satisfy
\begin{equation}
\label{U_j_x}
\left\{
\begin{array}{ll}
\d (U_{j+1})_x
&\d =pL_a'(|U_j+\e u_t^0|^{p-2}(U_j+\e u_t^0)(U_j+\e u_t^0)_x)\\
&\d \quad-(1+a)L_{a+1,b}'\left(|U_j+\e u_t^0|^p\frac{x}{\LR{x}}\frac{t+\LR{x}}{\LR{t+\LR{x}}}\right)\\
&\d \quad+(1+b)L_{a,b+1}'\left(|U_j+\e u_t^0|^p\frac{x}{\LR{x}}\frac{t-\LR{x}}{\LR{t-\LR{x}}}\right),\\
\d (U_1)_x&\d=0,
\end{array}
\right.
\end{equation}
so that the function space in which $\{U_j\}$ will converge is
\[
X:=\{U\in C^1(\R\times[0,T]) :\ \mbox{supp}\ U\subset\{|x|\le t+R\}\},
\]
equipping the norm
\[
\|U\|_X:=\|U\|+\|U_x\|,\quad \|U\|:=\sup_{\R\times[0,T]}|w(x,t)U(x,t)|, 
\]
where the weight function
\begin{equation}
\label{weight}
w(x,t):=\left\{
\begin{array}{llll}
\chi (1+t-|x|)^{1+a} + (1-\chi)&\mbox{if}& \ a>0,\\
\chi (1+t-|x|)^{1+a} + (1-\chi)\log^{-1}(t+|x|+R) &\mbox{if}& \ a=0,\\
\chi (1+t-|x|)^{1+a} + (1-\chi)(t+|x|+R)^{a} &\mbox{if}& \ -1 \le a<0,\\
\chi (1+t+|x|)^{1+a} + (1-\chi)(t+|x|+R)^{a} &\mbox{if}&\ a<-1.\\
\end{array}
\right.
\end{equation}
and the indicator function 
\begin{equation}
\label{chisubint}
\chi := \chi (x,t) :=
\left\{
\begin{array}{lll}
1 & t-|x| > R, \\
0 & \mbox{others}.
\end{array}
\right.
\end{equation}
First we note that $U_j\in X$ implies $U_{j+1}\in X$, namely $\{U_j\}$ is a sequence in $X$.
It is easy to check this fact by assumption on the initial data (\ref{supp_initial})
and the definitions of $L'_{a,b}$ in (\ref{integral_derivative}).

Then we have a following a priori estimates.
\begin{lem}
\label{lem:apriori1}
Let $U\in C(\R\times[0,T])$ and supp\ $U\subset\{(x,t)\in\R\times[0,T]:|x|\le t+R\}$. Then there exists a positive constant $C$ independent of $T$ and $\e$ such that
\begin{equation}
\label{apriori}
\|L'_{a,b}(|U|^p)\|\le CE_{a,b}(T)\|U\|^p,
\end{equation}
where
\begin{equation}
\label{E}
E_{a,b}(T):=
\left\{
\begin{array}{lllll}
1 &\mbox{for}\ a>0 \ \mbox{and} \ p(1+a)+b>0, \\
\log^p(T+3R) &\mbox{for}\  a=0\ \mbox{and}\ b\ge-p, \\
\log(T+3R) &\mbox{for}\  a>0\ \mbox{and}\ p(1+a)+b=0, \\
(T+2R)^{-ap} &\mbox{for}\ a<0 \ \mbox{and}\ b \ge -p,\\
(T+2R)^{-p(1+a)-b} &\mbox{for}\ p(1+a)+b<0\ \mbox{and}\ b<-p.\\
\end{array}
\right.
\end{equation}
\end{lem}
\begin{lem}
\label{lem:apriori2}
Suppose that the assumptions of Theorem \ref{thm:lower-bound} are fulfilled.
Assume that $U, U^0\in C(\R\times[0,T])$ with ${\rm supp}\ U\subset \{(x,t)\in \R \times[0,T]: \ |x|\le t+R\}$ and
 ${\rm supp}\ U^0\subset \{(x,t)\in \R \times[0,T]: \ (t-R)_{+}\le |x|\le t+R\}$ hold. 
Then, there exists a positive constant $C$ independent of $T$ and $\e$ such that
\begin{equation}
\label{apriori1}
\|L'_{a,b}(|U^0|^{p-j}|U|^j)\|\le CM(D(T)\|U\|)^j \quad \mbox{for}\ j=0,1,
\end{equation}
where $D_a(T)$ is defined by
\[
D_a(T):=
\left\{
\begin{array}{lll}
1 & \mbox{if}\ a>0,\\
\log(T+3R) & \mbox{if}\ a=0,\\
(T+2R)^{-a} & \mbox{if}\ a<0.
\end{array}
\right.
\]
\end{lem}
The proof of Lemma \ref{lem:apriori1} and Lemma \ref{lem:apriori2} are established in the next section.
Set
\[
M:=\|f'\|_{L^p(\R)}+\|f''\|_{L^p(\R)}+\|g\|_{L^p(\R)}+\|g'\|_{L^p(\R)}.
\]

\vskip10pt
\par\noindent
{\bf The convergence of the sequence $\v{\{U_j\}}$.}
\par
First we note that $\|U_1\|\le CM\e^p$ by Lemma \ref{lem:apriori1}.
Since (\ref{U_j}) and (\ref{apriori}) yield that
\[
\begin{array}{ll}
\|U_{j+1}\|
&\le \|L_{a,b}'(|U_j+\e u_t^0|^p)\|\\
&\le 2^{p-1}CM\e^p+2^{p-1}CE_{a,b}(T)\|U_j\|^p,
\end{array}
\]
the boundedness of $\{U_j\}$;
\begin{equation}
\label{bound_U}
\|U_j\|\le 2^pCM\e^p\quad(j\in\N)
\end{equation}
follows from
\begin{equation}
\label{condi1}
E_{a,b}(T)(2^pCM\e^p)^p\le M\e^p.
\end{equation}
Assuming (\ref{condi1}), one can estimate $U_{j+1}-U_j$ as follows.
\[
\begin{array}{ll}
\|U_{j+1}-U_j\| \hspace{-2.5mm}
&\le\|L'_{a,b}(|U_j+\e u_t^0|^p-|U_{j-1}+\e u_t^0|^p)\|\\
&\le 3^{p-1}p\|L_{a,b}'\left((|U_j|^{p-1}+|U_{j-1}|^{p-1}+\e^{p-1}|u_t^0|^{p-1})|U_j-U_{j-1}|\right)\|\\
&\le 3^{p-1}pCE_{a,b}(T)(\|U_j\|^{p-1}+\|U_{j-1}\|^{p-1})\|U_j-U_{j-1}\|\\
&\quad + 3^{p-1}pCM\e^{p-1}D_a(T)\|U_j-U_{j-1}\| \\
&\le 3^{p-1}pCE_{a,b}(T)2(2^p3CM\e^p)^{p-1}\|U_j-U_{j-1}\| \\
&\quad + 3^{p-1}pCM\e^{p-1}D_a(T)\|U_j-U_{j-1}\| 
\end{array}
\]
Therefore the convergence of $\{U_j\}$ follows from
\begin{equation}
\label{convergence}
\|U_{j+1}-U_j\|\le\frac{1}{2}\|U_j-U_{j-1}\|\quad(j\ge2)
\end{equation}
provided (\ref{condi1}) and
\begin{equation}
\label{condi2}
3^{p-1}pCE_{a,b}(T)2(2^pCM\e^p)^{p-1}\le\frac{1}{4}
\end{equation}
and
\begin{equation}
\label{condi3}
3^{p-1}pCM\e^{p-1}D_a(T)\le\frac{1}{4}
\end{equation}
are fulfilled.

\vskip10pt
\par\noindent
{\bf The convergence of the sequence $\v{\{(U_j)_x\}}$.}
\par
First we note that $\|(U_1)_x\|\le CM\e^p$ by Lemma \ref{lem:apriori1}.
Assume that (\ref{condi1}), (\ref{condi2}) and (\ref{condi3}) are fulfilled.
It follows from (\ref{U_j_x}) and (\ref{apriori}) that
\[
\begin{array}{ll}
\d \|(U_{j+1})_x\|
&\d \le \|L_{a,b}'\left(|U_j+\e u_t^0|^{p-1}|(U_j+\e u_t^0)_x|\right)\| \\
&\d \quad +|1+a|\|L'_{a+1,b}(|U_j+\e u_t^0|^p)\| \\
&\d \quad +|1+b|\|L'_{a,b+1}(|U_j+\e u_t^0|^p)\| \\
&\d \le 2^{p-1}CE_{a,b}(T)\|U_j\|^{p-1}\|(U_j)_x\| +2^{p-1}CM\e^p\\
&\d \quad+|1+a|2^{p-1}CE_{a+1,b}(T)\|U_j\|^p +2^{p-1}CM\e^p\\
&\d \quad+|1+b|2^{p-1}CE_{a,b+1}(T)\|U_j\|^p +2^{p-1}CM\e^p\\
&\d \le 2^{p-1}CE_{a,b}(T)(2^p3CM\e^p)^{p-1}\|(U_j)_x\| + 2^{p-1}3CM\e^p\\
&\d \quad+\{|1+a|2^{p-1}CE_{a+1,b}(T)+|1+b|2^{p-1}CE_{a,b+1}(T)\}(2^p3CM\e^p)^p.
\end{array}
\]
Hence the boundedness of $\{(U_j)_x\}$;
\begin{equation}
\label{bound_U_x}
\|(U_j)_x\|\le 2^p3CM\e^p\quad(j\in\N)
\end{equation}
follows from
\begin{equation}
\label{boundedness_condi}
\{E_{a,b}(T)+|1+a|E_{a+1,b}(T)+|1+b|E_{a,b+1}(T)\}(2^p3CM\e^p)^p \le 3M\e^p
\end{equation}
Assuming (\ref{boundedness_condi}), one can estimate $(U_{j+1})_x-(U_j)_x$ as follows.
\[
\begin{array}{ll}
\|(U_{j+1})_x-(U_j)_x\|
&\le\|L'_{a,b}(|U_j+\e u_t^0|^{p-2}(U_j+\e u_t^0)(U_j+\e u_t^0)_x \\ 
& \qquad -|U_{j-1}+\e u_t^0|^{p-2}(U_{j-1}+\e u_t^0)(U_{j-1}+\e u_t^0)_x)\|\\
&\quad+|1+a|\|L'_{a+1,b}\left(|U_j+\e u_t^0|^p-|U_{j-1}+\e u_t^0|^p\right)\| \\
&\quad+|1+b|\|L'_{a,b+1}\left(|U_j+\e u_t^0|^p-|U_{j-1}+\e u_t^0|^p\right)\|.
\end{array}
\]
The first term on the right hand side of above inequality is split into two pieces
according to
\[
\begin{array}{ll}
|U|^{p-2}UU_x-|U'|^{p-2}U'U'_x \\
=(|U|^{p-2}U-|U'|^{p-2}U')U_x+|U'|^{p-2}U'(U_x-U'_x).
\end{array}
\]
Since
\[
\begin{array}{l}
\left||U_j+\e u_t^0|^{p-2}(U_j+\e u_t^0)-|U_{j-1}+\e u_t^0|^{p-2}(U_{j-1}+\e u_t^0)\right| 
\\
\le
\left\{
\begin{array}{ll}
(p-1)(|U_j+\e u_t^0|^{p-2}+|U_{j-1}+\e u_t^0|^{p-2})|U_j-U_{j-1}| & \mbox{when}\ p\ge2,\\
2|U_j-U_{j-1}|^{p-1} & \mbox{when}\ 1<p<2,\\
\end{array}
\right.
\end{array}
\]
and
\[
\begin{array}{ll}
(|U_j+\e u_t^0|^{p-2}+|U_{j-1}+\e u_t^0|^{p-2})|(U_j+\e u_t^0)_x| \\
\le 2^p(|U_j|^{p-1}+|U_{j-1}|^{p-1}+|(U_j)_x|^{p-1}+\e^{p-1}|u_t^0|^{p-1}+\e^{p-1}|u_{tx}^0|^{p-1}),
\end{array}
\]
the similar manner of  handling $L_{a+1,b}'$ and $L_{a,b+1}'$ to above computations leads to
\[
\begin{array}{l}
\|(U_{j+1})_x-(U_j)_x\|\\
\le 
\left\{
\begin{array}{ll}
F \ \|U_j-U_{j-1}\| & \mbox{when}\ p\ge2,\\
2 \|U_j-U_{j-1}\|^{p-1} & \mbox{when}\ 1<p<2,\\
\end{array}
\right.
\\
\quad+2^pC\{E_{a,b}(T)\|U_{j-1}\|^{p-1}+M\e^{p-1}D_a(T)\}\|(U_j)_x-(U_{j-1})_x\|\\
\quad+|1+a|2^pCp\{E_{a+1,b}(T)\|U_{j-1}\|^{p-1}+M\e^{p-1}D_{a+1}(T)\}\|U_j-U_{j-1}\| \\
\quad+|1+b|2^pCp\{E_{a,b+1}(T)\|U_{j-1}\|^{p-1}+M\e^{p-1}D_a(T)\}\|U_j-U_{j-1}\|,
\end{array}
\]
where 
\[
F := 2^p(p-1) \{ CE_{a,b}(T)(\|U_j\|^{p-1}+\|U_{j-1}\|^{p-1}+\|(U_j)_x\|^{p-1}) +2CM\e^{p-1}D_a(T)\}
\]
Hence it follows from (\ref{convergence}) that
\[
\begin{array}{ll}
\|(U_{j+1})_x-(U_j)_x\|
&\le 2^pC\{E_{a,b}(T)(2^p3CM\e^p)^{p-1}+M\e^{p-1}D_a(T)\} \times \\
&\quad \times \|(U_j)_x-(U_{j-1})_x\|\\
&\quad\d+
 O\left(\frac{1}{2^{j \min \{ (p-1),1\}} }\right)
\end{array}
\]
as $j\rightarrow\infty$.
Here we have employed the fact that $E_{a+1,b}(T)$ and $E_{a,b+1}(T)$ 
is dominated by $E_{a,b}(T)$ with some positive constant, and that $D_{a+1}(T)$ is dominated by $D_a(T)$.
Therefore we obtain the convergence of $\{(U_j)_x\}$ provided
\begin{equation}
\label{condi4}
2^{p-1}C\{E_{a,b}(T)(2^p3CM\e^p)^{p-1}+2M\e^{p-1}D_a(T)\}\le\frac{1}{2}.
\end{equation}
%
%
%
\hfill$\Box$


\section{Proof of Lemma \ref{lem:apriori1} and Lemma \ref{lem:apriori2}}
\par
In this section, we prove priori estimates (\ref{apriori1}) and (\ref{apriori}).
Recall the definition of $L'_{a,b}$ in (\ref{nonlinear_derivative}).
From now on, a positive constant $C$ independent of $T$ and $\e$
may change from line to line.
For any $(x,t) \in \R \times [0,\infty)$, we have
\[
3^{-1}(1+|t-|x||) \leq \LR{t-\LR{x}} \leq \sqrt{2}  (1+|t-|x||) 
\]
and
\[
2^{-1}(1+t+|x|) \leq  \LR{t+\LR{x}} \leq \sqrt{2} (1+t+|x|).
\]
First, Lemma \ref{lem:apriori1} is shown. 
Let us assume the conditions of Lemma \ref{lem:apriori1}

Therefore, we have that
\[
|L_{a,b}'(|U|^p)(x,t)|\le
C\|U\|^p\{I_+(x,t)+I_-(x,t)\},
\]
where the integrals $I_+$ and $I_-$ are defined by
\[
I_\pm(x,t):=\int_0^t\frac{w^{-p}(x\pm(t-s) , s)\chi_\pm(x,t;s)}{(1+s+|t-s\pm x|)^{1+a}(1+|s-|t-s\pm x||)^{1+b}}ds
\]
and the characteristic functions $\chi_+$ and $\chi_-$ are defined by
\[
\begin{array}{ll}
\chi_\pm(x,t;s)
&:=\chi_{\{s: |t-s\pm x|\le s+R\}}\\
&=
\left\{
\begin{array}{ll}
1 & \mbox{when $s$ satisfies }|t-s\pm x|\le s+R,\\
0 & \mbox{otherwise},
\end{array}
\right.
\end{array}
\]
respectively.
First we note that it is sufficient to  estimate $I_\pm$ for $x\ge0$ due to its symmetry,
\[
I_+(-x,t)=I_-(x,t).
\]
Hence it follows from $0\le x\le t+R$ as well as
\[
|t-s+x|\le s+R\quad\mbox{and}\quad 0\le s\le t
\]
that
\[
\max\left\{0, \frac{t+x-R}{2}\right\}\le s\le t.
\]
In the case of $t+x<R$, since $t \le R$ and $x \le R$, $I_\pm$ is suppressed by a constant from above.
In the case of $t+x \ge R$,
\[
\begin{array}{ll}
\d I_+(x,t)
&\d = \int_{(t+x-R)/2}^t\frac{w^{-p}(x+t-s , s)ds}{(1+s+|t-s+x|)^{1+a}(1+|s-|t-s+x||)^{1+b}} \\
&\d = \int_{(t+x-R)/2}^t\frac{w^{-p}(x+t-s , s)ds}{(1+t+x)^{1+a}(1+|2s-t-x|)^{1+b}}
\end{array}
\]
holds. If $(t-R)_{+}\le x \le t+R$, it follows from 
\[
\begin{array}{ll}
& \d \int_{(t+x-R)/2}^t \frac{ds}{(1+|2s-t-x|)^{1+b}} \\
& \d \le \int_{(t+x-R)/2}^{(t+x+R)/2} \frac{ds}{(1-|2s-t-x|)^{1+b}} \le R \max\{1,(1+R)^{-1-b}\} \\
\end{array}
\]
that
\[
\begin{array}{ll}
 I_+(x,t) 
& \d = \int_{(t+x-R)/2}^t
\left\{
\begin{array}{ll}
1 & (a>0)\\
\log^p(s+|t-s+x|+R) & (a=0) \\
(s+|t-s+x|+R)^{-ap} & (a<0)
\end{array}
\right\} 
\times \\
& \quad  \times 
\d \frac{ds}{(1+t+x)^{1+a}(1+|2s-t-x|)^{1+b}} \\
&\d \le \frac{C}{(1+t+x)^{1+a}}
\left\{
\begin{array}{ll}
1 & (a>0)\\
\log^p(t+x+R) & (a=0) \\
(t+x+R)^{-ap} & (a<0)
\end{array}
\right. \\
&\d \le \frac{Cw^{-1}(x,t)}{(1+t+x)^{1+a}}
\left\{
\begin{array}{ll}
1 & (a>0)\\
\log^{p-1}(t+x+R) & (a=0) \\
(t+x+R)^{-a(p-1)} & (a<0)
\end{array}
\right. \\
&\d \le CE_{a,b}(T) w^{-1}(x,t).
\end{array}
\]
Therefore, we can assume that $t \ge x+R$. 
In this case, $I_+$ can be divided as follows; 
\[
\begin{array}{ll}
\d I_+(x,t)
&\d = \int_{(t+x+R)/2}^t\frac{w^{-p}(x+t-s , s)ds}{(1+t+x)^{1+a}(1+|2s-t-x|)^{1+b}} \\
& \d \quad + \int_{(t+x-R)/2}^{(t+x+R)/2}\frac{w^{-p}(x+t-s , s)ds}{(1+t+x)^{1+a}(1+|2s-t-x|)^{1+b}}. \\
\end{array}
\]
When $a>0$, we have
\[
\begin{array}{ll}
\d I_+(x,t)
&\d = \int_{(t+x+R)/2}^t\frac{(1+s-|t-s+x|)^{-p(1+a)}ds}{(1+t+x)^{1+a}(1+|2s-t-x|)^{1+b}} \\
& \d \quad + \int_{(t+x-R)/2}^{(t+x+R)/2}\frac{ds}{(1+t+x)^{1+a}(1+|2s-t-x|)^{1+b}}\\
&\d \le C \int_{(t+x+R)/2}^t\frac{ds}{(1+t+x)^{1+a}(1+2s-t-x)^{p(1+a)+1+b}} \\
& \d \quad + \int_{(t+x-R)/2}^{(t+x+R)/2}\frac{ds}{(1+t+x)^{1+a}(1+|2s-t-x|)^{1+b}}\\
&\d \le \frac{C}{(1+t+x)^{1+a}}
\left\{
\begin{array}{lll}
1 & \mbox{for} \ p(1+a)+b > 0\\
\log(1+t-x) & \mbox{for} \ p(1+a)+b = 0\\
(1+t-x)^{-p(1+a)-b} & \mbox{for} \ p(1+a)+b < 0\\
\end{array} 
\right.\\
& \d \quad + \frac{C}{(1+t+x)^{1+a}}.
\end{array}
\]
When $a=0$, we have
\[
\begin{array}{ll}
\d I_+(x,t)
&\d = \int_{(t+x+R)/2}^t\frac{(1+s-|t-s+x|)^{-p(1+a)}ds}{(1+t+x)^{1+a}(1+|2s-t-x|)^{1+b}} \\
& \d \quad + \int_{(t+x-R)/2}^{(t+x+R)/2}\frac{\log^p(s+|t-s+x|+R)ds}{(1+t+x)^{1+a}(1+|2s-t-x|)^{1+b}}\\
&\d \le \frac{C}{1+t+x}
\left\{
\begin{array}{lll}
1 & \mbox{for} \ p+b > 0\\
\log(1+t-x) & \mbox{for} \ p+b = 0\\
(1+t-x)^{-p-b} & \mbox{for} \ p+b < 0\\
\end{array} 
\right.\\
& \d \quad + \frac{C\log^p(t+x+R)}{1+t+x}.
\end{array}
\]
When $-1\le a<0$, we have
\[
\begin{array}{ll}
\d I_+(x,t)
&\d = \int_{(t+x+R)/2}^t\frac{(1+s-|t-s+x|)^{-p(1+a)}ds}{(1+t+x)^{1+a}(1+|2s-t-x|)^{1+b}} \\
& \d \quad + \int_{(t+x-R)/2}^{(t+x+R)/2}\frac{(s+|t-s+x|+R)^{-ap}ds}{(1+t+x)^{1+a}(1+|2s-t-x|)^{1+b}}\\
&\d \le \frac{C}{(1+t+x)^{1+a}}
\left\{
\begin{array}{lll}
1 & \mbox{for} \ p(1+a)+b > 0\\
\log(1+t-x) & \mbox{for} \ p(1+a)+b = 0\\
(1+t-x)^{-p(1+a)-b} & \mbox{for} \ p(1+a)+b < 0\\
\end{array} 
\right.\\
& \d \quad + \frac{C(t+x+R)^{-ap}}{(1+t+x)^{1+a}}.
\end{array}
\]
When $a<-1$, we have
\[
\begin{array}{ll}
\d I_+(x,t)
&\d = \int_{(t+x+R)/2}^t\frac{(1+s+|t-s+x|)^{-p(1+a)}ds}{(1+t+x)^{1+a}(1+|2s-t-x|)^{1+b}} \\
& \d \quad + \int_{(t+x-R)/2}^{(t+x+R)/2}\frac{(s+|t-s+x|+R)^{-ap}ds}{(1+t+x)^{1+a}(1+|2s-t-x|)^{1+b}}\\
&\d \le \frac{C}{(1+t+x)^{(p+1)(1+a)}}
\left\{
\begin{array}{lll}
1 & \mbox{for} \ b > 0\\
\log(1+t-x) & \mbox{for} \ b = 0\\
(1+t-x)^{-b} & \mbox{for} \ b < 0\\
\end{array} 
\right.\\
& \d \quad + \frac{C(t+x+R)^{-ap}}{(1+t+x)^{1+a}}.
\end{array}
\]
The following inequality can be obtained by paying attention to the size of the powers:
\[
I_+(x,t) \le CE_{a,b}(T) w^{-1}(x,t).
\]
\par
On the other hand, the estimate for $I_-$ is divided into two cases.
If $t-x\ge R$, then  $|t-s-x|\le s+R$ yields that
\[
\begin{array}{ll}
I_-(x,t)
& \d \le\int_{(t-x-R)/2}^{t-x}\frac{w^{-p}(x-t+s , s)ds}{(1+t-x)^{1+a}(1+|2s-t+x|)^{1+b}}ds \\
& \d \quad +\int_{t-x}^t\frac{w^{-p}(x-t+s , s)ds}{(1+|2s-t+x|)^{1+a}(1+t-x)^{1+b}}ds
\end{array} 
\]
follows. Furthermore, the first integral is separated by the domain of the weight function $w$.
When $a>0$, we have
\[
\begin{array}{ll}
I_-(x,t)
& \d \le\int_{(t-x-R)/2}^{(t-x+R)/2}\frac{ds}{(1+t-x)^{1+a}(1+|2s-t+x|)^{1+b}}ds \\
& \d \quad +\int_{(t-x+R)/2}^{t-x}\frac{(1+2s-t+x)^{-p(1+a)}ds}{(1+t-x)^{1+a}(1+2s-t+x)^{1+b}}ds \\
& \d \quad +\int_{t-x}^t\frac{(1+t-x)^{-p(1+a)}ds}{(1+2s-t+x)^{1+a}(1+t-x)^{1+b}}ds \\
& \d \le \frac{C}{(1+t-x)^{1+a}} \\
& \d \quad +\frac{C}{(1+t-x)^{1+a}} 
\left\{
\begin{array}{lll}
1 & \mbox{for} \ p(1+a)+b > 0\\
\log(1+t-x) & \mbox{for} \ p(1+a)+b  = 0\\
(1+t-x)^{-p(1+a)-b} & \mbox{for} \ p(1+a)+b  < 0\\
\end{array}
\right. \\
& \d \quad + \frac{C}{(1+t-x)^{(p+1)(1+a)+b}}  \\
\end{array} 
\]
When $a=0$, we have
\[
\begin{array}{ll}
I_-(x,t)
& \d \le\int_{(t-x-R)/2}^{(t-x+R)/2}\frac{\log^{p}(t-x+R)ds}{(1+t-x)^{1+a}(1+|2s-t+x|)^{1+b}}ds \\
& \d \quad +\int_{(t-x+R)/2}^{t-x}\frac{(1+2s-t+x)^{-p(1+a)}ds}{(1+t-x)^{1+a}(1+|2s-t+x|)^{1+b}}ds \\
& \d \quad +\int_{t-x}^t\frac{(1+t-x)^{-p(1+a)}ds}{(1+|2s-t+x|)^{1+a}(1+t-x)^{1+b}}ds \\
\end{array}
\]
Regarding the above, for the first integral,
\[
\int_{(t-x-R)/2}^{(t-x+R)/2}\frac{\log^{p}(t-x+R)ds}{(1+t-x)(1+|2s-t+x|)^{1+b}}ds \le \frac{C}{1+t-x} \log^{p}(t-x+R)
\]
hold. For the second intedral, 
\[
\begin{array}{ll}
& \d \int_{(t-x+R)/2}^{t-x}\frac{(1+2s-t+x)^{-p}ds}{(1+t-x)(1+|2s-t+x|)^{1+b}}ds \\
& \d \quad \le \frac{C}{1+t-x} 
\left\{
\begin{array}{lll}
1 & \mbox{for} \ p+b > 0\\
\log(1+t-x) & \mbox{for} \ p+b  = 0\\
(1+t-x)^{-p-b} & \mbox{for} \ p+b  < 0\\
\end{array}
\right. 
\end{array}
\]
hold. For the third integral, by the inequality; 
\[
(1+t-x)^{-p-b}\log\left( \frac{1+t+x}{1+t-x} \right) \le C (1+t+x)^{-p-b} \quad \mbox{for} \ p+b<0,
\]
it follows that 
\[
\begin{array}{ll}
& \d \int_{t-x}^t\frac{(1+t-x)^{-p}ds}{(1+|2s-t+x|)(1+t-x)^{1+b}}ds  \\
& \d \quad \le \frac{C}{(1+t-x)^{1+p+b}}\log\left( \frac{1+t+x}{1+t-x} \right) \\
& \d \quad \le \frac{C}{1+t-x} 
\left\{
\begin{array}{lll}
\log(1+t+x) & \mbox{for} \ p+b  \ge 0,\\
(1+t+x)^{-p-b} & \mbox{for} \ p+b  < 0.\\
\end{array}
\right. 
\end{array}
\]
From the above, we have 
\[
I_-(x,t) \le 
\frac{C}{1+t-x} 
\left\{
\begin{array}{lll}
\log^p(t+x+R) & \mbox{for} \ p+b  \ge 0,\\
(1+t+x)^{-p-b} & \mbox{for} \ p+b  < 0.\\
\end{array}
\right. 
\]
When $-1<a<0$, we have
\[
\begin{array}{ll}
I_-(x,t)
& \d \le\int_{(t-x-R)/2}^{(t-x+R)/2}\frac{(t-x+R)^{-ap}ds}{(1+t-x)^{1+a}(1+|2s-t+x|)^{1+b}}ds \\
& \d \quad +\int_{(t-x+R)/2}^{t-x}\frac{(1+2s-t+x)^{-p(1+a)}ds}{(1+t-x)^{1+a}(1+|2s-t+x|)^{1+b}}ds \\
& \d \quad +\int_{t-x}^t\frac{(1+t-x)^{-p(1+a)}ds}{(1+|2s-t+x|)^{1+a}(1+t-x)^{1+b}}ds \\
& \d \le \frac{C(t-x+R)^{-ap}}{(1+t-x)^{1+a}} \\
& \d \quad +\frac{C}{(1+t-x)^{1+a}} 
\left\{
\begin{array}{lll}
1 & \mbox{for} \ p(1+a)+b > 0\\
\log(1+t-x) & \mbox{for} \ p(1+a)+b  = 0\\
(1+t-x)^{-p(1+a)-b} & \mbox{for} \ p(1+a)+b  < 0\\
\end{array}
\right. \\
& \d \quad + \frac{C(1+t+x)^{-a} }{(1+t-x)^{p(1+a)+1+b}}  \\
\end{array} 
\]
When $a<-1$, we have
\[
\begin{array}{ll}
I_-(x,t)
& \d \le\int_{(t-x-R)/2}^{(t-x+R)/2}\frac{(t-x+R)^{-ap}ds}{(1+t-x)^{1+a}(1+|2s-t+x|)^{1+b}}ds \\
& \d \quad +\int_{(t-x+R)/2}^{t-x}\frac{(1+t-x)^{-p(1+a)}ds}{(1+t-x)^{1+a}(1+|2s-t+x|)^{1+b}}ds \\
& \d \quad +\int_{t-x}^t\frac{(1+2s-t+x)^{-p(1+a)}ds}{(1+|2s-t+x|)^{1+a}(1+t-x)^{1+b}}ds \\
& \d \le \frac{C(t-x+R)^{-ap}}{(1+t-x)^{1+a}} \\
& \d \quad +\frac{C(1+t-x)^{-p(1+a)}}{(1+t-x)^{1+a}} 
\left\{
\begin{array}{lll}
1 & \mbox{for} \ b > 0\\
\log(1+t-x) & \mbox{for} \ b  = 0\\
(1+t-x)^{-b} & \mbox{for} \ b  < 0\\
\end{array}
\right. \\
& \d \quad + \frac{C}{(1+t-x)^{1+b}}(1+t+x)^{-(p+1)(1+a)+1}   \\
\end{array} 
\]
Therefore we obtain
\[
I_-\le CE_{a,b}(T)\quad\mbox{in}\ \R\times[0,T]\cap\{t-x\ge R \}.
\]
If $(-R\le)t-x\le R$, $|t-s-x|\le s+R$ yields that
\[
\begin{array}{ll}
I_-(x,t)
& \d \le\int_0^{(t-x)_+}\frac{w^{-p}(x-t+s , s)ds}{(1+|t-x|)^{1+a}(1+|2s-t+x|)^{1+b}}ds \\
& \d \quad +\int_{(t-x)_+}^t\frac{w^{-p}(x-t+s , s)ds}{(1+|2s-t+x|)^{1+a}(1+|t-x|)^{1+b}}ds.
\end{array} 
\]
Since $t-x$ is equivalent to $R$, the first term is evaluated with a constant. Therefore,  the classification can be based on the weight of the second term.
When $a<0$, we have
\[
I_-(x,t)\le C + C(1+x)^{-a}\le C(T+2R)^{-a}.
\]
When $a=0$, we have
\[
I_-(x,t)\le C + C\log\frac{1+x}{1-t+x}\le C \log(T+2R).
\]
When $a>0$, we have
\[
I_-(x,t)\le C + C(1-t+x)^{-a}\le C.
\]
Therefore we obtain
\[
I_-\le CE_{a,b}(T)\quad\mbox{in}\ \R\times[0,T]\cap\{-R\le t-x\le0\}.
\]
Summing up all the estimates for $I_+$ and $I_-$, we have
\[
|L_a'(|U|^p)|\le
C\|U\|^pE_a(T)\quad\mbox{in}\ \R\times[0,T].
\]
This completes the proof of Lemma \ref{lem:apriori1}.

Next is the proof for Lemma \ref{lem:apriori2}. For 
\[
|L_a'(|U^0|^{p-j}|U|^j)|\le
C\|U\|^j\{J_++J_-\}, 
\]
the proof can be obtained by the same calculation as in Lemma \ref{lem:apriori1}, assuming that $J_{\pm}$ is defined by changing the exponent of $I_{\pm}$ from $-p$ to $-j$ and from the function $\chi_\pm$ to
\[
\begin{array}{ll}
\chi^0_\pm(x,t;s)
&:=\chi_{\{s: (s-R)_+ \le |t-s\pm x|\le s+R\}}\\
&=
\left\{
\begin{array}{ll}
1 & \mbox{when $s$ satisfies }(s-R)_+ \le |t-s\pm x|\le s+R,\\
0 & \mbox{otherwise}.
\end{array}
\right.
\end{array}
\]
\hfill$\Box$.

\section{Proof of Theorem \ref{thm:upper-bound}}
\par
In this section, a positive constant $C$ independent of $T$ and $\e$
may change from line to line.
Let $U\in C(\R\times[0,T])$ with
\begin{equation}
\label{supp_integralsol}
\mbox{supp}\ U(x,t)\subset\{(x,t)\in\R\times[0,T]:|x|\le t+R\}
\end{equation}
 be a solution of the integral equation (\ref{integral_derivative}),
namely
\[
U=\e u_t^0+L'_{a,b}(|U|^p).
\]
Then it is easy to see by simple integration that
\[
V(x,t):=\int_0^tU(x,s)ds+\e f(x)
\]
satisfies a integral equation,
\[
V=\e u^0+L_a(|V_t|^p).
\]

\par
Set $t=x+R,\ x\ge R$. Then, inverting the order of the $(y,s)$-integral and
diminishing its domain, we have that
\begin{equation}
\label{b>0blowup}
\begin{array}{ll}
\d L_a(|V_t|^p)(x,t) 
& \d \ge C\int_R^xdy\int_{y-R}^{y+R}\frac{|V_t(y,s)|^p}{(1+s+y)^{1+a}(1+|s-y|)^{1+b}}ds \\
& \d \ge C\int_R^xdy\int_{y-R}^{y+R}\frac{|V_t(y,s)|^p}{(1+2y+R)^{1+a}(1+R)^{1+b}}ds \\
& \d \ge \frac{C}{2^{1+a}(1+R)^{1+b}}\int_R^xdy\frac{1}{(R+y)^{1+a}}\int_{y-R}^{y+R}|V_t(y,s)|^pds. \\
\end{array}
\end{equation}
The integral inequality is consistent with one of the solution explosion proof in \cite{KMT22} using the method in Zhou \cite{Zhou01} except for the difference in constants, so we omit the following calculations and obtain the following result.

\begin{prop}[\cite{KMT22}]
\label{prop:1}
Assume (\ref{supp_initial})
and
\begin{equation}
\label{positive_non-zero1}
\int_{\R}g(x)dx >0.
\end{equation}
Then, there exists a positive constant $\e_1=\e_1(g,p,a,R)>0$ such that
a solution $u_t\in C(\R\times[0,T])$ with
\[
\mbox{\rm supp}\ u_t(x,t)\subset\{(x,t)\in\R\times[0,T]:|x|\le t+R\}
\]
of associated integral equations (\ref{integral_derivative})
to (\ref{IVPderivative})
cannot exist whenever $T$ satisfies
\begin{equation}
\label{upper-bound_a=0}
T\ge
\left\{
\begin{array}{ll}
C\e^{-(p-1)/(-a)} & \mbox{for}\ a<0,\\
\exp\left(C\e^{-(p-1)}\right) & \mbox{for}\ a=0,
\end{array}
\right.
\end{equation}
where $0<\e\le\e_1$, $C$ is a positive constant independent of $\e$.
\end{prop}
Therefore, in order to prove the Theorem \ref{thm:upper-bound}, it suffices to prove the following proposition by point-wise estimates of the solution.
\begin{prop}
\label{prop:2}
Assume (\ref{supp_initial})
and
\begin{equation}
\label{positive_non-zero2}
f(x) \equiv 0 , \quad  \int_{\R}g(x)dx >0.
\end{equation}
Then, there exists a positive constant $\e_2=\e_2(g,p,a,b,R)>0$ such that
a solution $u_t\in C(\R\times[0,T])$ with
\[
\mbox{\rm supp}\ u_t(x,t)\subset\{(x,t)\in\R\times[0,T]:|x|\le t+R\}
\]
of associated integral equations (\ref{integral_derivative})
to (\ref{IVPderivative})
cannot exist whenever $T$ satisfies
\begin{equation}
\label{upper-bound_p(1+a)+b=0}
T\ge
\left\{
\begin{array}{ll}
C\e^{-p(p-1)/(-p(1+a)+b)} & \mbox{for}\ p(1+a)+b<0,\\
\exp\left(C\e^{-p(p-1)}\right) & \mbox{for}\ p(1+a)+b=0,
\end{array}
\right.
\end{equation}
where $0<\e\le\e_2$, $C$ is a positive constant independent of $\e$.
\end{prop}

We will prove the Proposition \ref{prop:2} below. 
Assume that $U$ is a solution of associated integral eqations (\ref{integral_derivative}) with 
\[
\mbox{\rm supp}\ U(x,t)\subset\{(x,t)\in\R\times[0,T]:|x|\le t+R\}
\]
and $(x,t)\in D$, where
\begin{equation}
\label{D}
D:=\{(x,t)\in \R \times[0,T]: t+|x|\ge R,\   t-|x|\ge R\}.
\end{equation}
From associated integral eqations (\ref{integral_derivative}), we have 
\[
\begin{array}{lll}
\d U(x,t) 
& \d \ge  \e u^0_t(x,t) \\
& \d \quad +C_0\Big\{\int_{(t-x-R)/2}^{(t-x+R)/2} &\d \frac{|U(x-t+s,s)|^p}{(1+t-x)^{1+a}(1+|2s-t+x|)^{1+b}}ds \\
& \d \qquad \quad +\int_{(t-x+R)/2}^{t-x} &\d \frac{|U(x-t+s,s)|^p}{(1+t-x)^{1+a}(1+|2s-t+x|)^{1+b}}ds \\
& \d \qquad \quad+\int_{t-x}^t &\d \frac{|U(x-t+s,s)|^p}{(1+|2s-t+x|)^{1+a}(1+t-x)^{1+b}}ds \\
& \d \qquad \quad+\int_{(t+x+R)/2}^t &\d \frac{|U(x+t-s,s)|^p}{(1+t+x)^{1+a}(1+|2s-t-x|)^{1+b}}ds \\
& \d \qquad \quad+ \int_{(t+x-R)/2}^{(t+x+R)/2} &\d \frac{|U(x+t-s,s)|^p}{(1+t+x)^{1+a}(1+|2s-t-x|)^{1+b}}ds \Big\},  \\
\end{array}
\]
where
\[
C_0 := \frac{3^{-|a|}2^{-|b|}}{8\sqrt{2}}.
\]
Since $U(x,t) \ge \e u^0_t(x,t) = 2^{-1}\e \{g(x+t)+g(x-t)\}$, 
\[
\begin{array}{lll}
& \d \int_{(t-x-R)/2}^{(t-x+R)/2} \frac{|U(x-t+s,s)|^p}{(1+t-x)^{1+a}(1+|2s-t+x|)^{1+b}}ds \\
& \d \ge \int_{(t-x-R)/2}^{(t-x+R)/2} \frac{2^{-p}\e^p \{g(x-t+2s)+g(x-t)\}^p}{(1+t-x)^{1+a}(1+|2s-t+x|)^{1+b}}ds \\
& \d \ge \frac{C\|g\|_{L^p}^p \e^p}{2^p(1+t-x)^{1+a}} =: \frac{C_g \e^p}{(1+t-x)^{1+a}} 
\end{array}
\]
holds. Therefore, the following integral inequality is obtained for use in the iterative method of each point evaluation: 
\begin{equation}\label{ineq1}
\begin{array}{lll}
\d U(x,t) 
& \d\ge  C_0 \int_{(t-x+R)/2}^{t-x}  \frac{|U(x-t+s,s)|^p}{(1+t-x)^{1+a}(1+|2s-t+x|)^{1+b}}ds \\
& \d  \quad+C_0\int_{(t+x+R)/2}^t  \frac{|U(x+t-s,s)|^p}{(1+t+x)^{1+a}(1+|2s-t-x|)^{1+b}}ds \\
& \d \quad +\frac{C_g \e^p}{(1+t-x)^{1+a}}
\end{array}
\end{equation}
\par\noindent
{\bf Case 1. $\v{p(1+a)+b=0}$.} 
\par\noindent
Assume that an estimate
\begin{equation}
\label{case_p(1+a)+b=0}
U(x,t)\ge \frac{M_n}{(1+t-x)^{1+a}}\left\{\log\left(\frac{1+t-x}{1+R}\right)\right\}^{a_n}
\quad\mbox{in}\ D
\end{equation}
holds, where $a_n\ge0$ and $M_n>0$.
The sequences $\{a_n\}$ and $\{M_n\}$ are defined later. 
Then it follows from (\ref{ineq1}) of the first integral term and (\ref{case_p(1+a)+b=0}) that
\[
\begin{array}{lll}
\d U(x,t) 
& \d \ge C_0\int_{(t+x+R)/2}^t  \frac{|U(x+t-s,s)|^p}{(1+t+x)^{1+a}(1+|2s-t-x|)^{1+b}}ds \\
& \d \ge C_0M_n^p\int_{(t+x+R)/2}^t  \frac{1}{(1+t+x)^{1+a}(1+2s-t-x)^{p(1+a)+1+b}} \times \\
& \d \quad \times \left\{\log\left(\frac{1+2s-t-x}{1+R}\right)\right\}^{pa_n}ds \\
& \d = \frac{C_0M_n^p}{(pa_n+1)(1+t+x)^{1+a}}\left\{\log\left(\frac{1+t-x}{1+R}\right)\right\}^{pa_n+1}.
\end{array}
\]
Next, from (\ref{ineq1}) of the second integral term, 
\[
\begin{array}{lll}
\d U(x,t) 
& \d\ge  C_0 \int_{(t-x+R)/2}^{t-x}  \frac{|U(x-t+s,s)|^p}{(1+t-x)^{1+a}(1+|2s-t+x|)^{1+b}}ds \\
& \d\ge  \frac{C_0^{p+1}M_n^{p^2}}{(pa_n+1)^p} \int_{(t-x+R)/2}^{t-x}  \frac{1}{(1+t-x)^{1+a}(1+2s-t+x)^{p(1+a)+1+b}}\times  \\
& \d \quad \times \left\{\log\left(\frac{1+t-x}{1+R}\right)\right\}^{p^2a_n+p}ds \\
& \d = \frac{C_0^{p+1}M_n^{p^2}}{(pa_n+1)^p(1+t-x)^{1+a}}\left\{\log\left(\frac{1+t-x}{1+R}\right)\right\}^{p^2a_n+p+1} \\
\end{array}
\]
holds. Therefore, if $\{a_n\}$ is defined by
\begin{equation}
\label{a_n}
a_{n+1}=p^2a_n+p+1,\ a_1=0,
\end{equation}
then (\ref{case_p(1+a)+b=0}) holds for all $n\in\N$ as far as $M_n$ satisfies
\begin{equation}
\label{M_n}
M_{n+1}\le C_0M_n^{p^2}.
\end{equation}
In view of (\ref{ineq1}), we note that (\ref{case_p(1+a)+b=0}) holds for $n=1$ with
\begin{equation}
\label{M_1}
M_1=C_g\e^p.
\end{equation}
Therefore, it follows from (\ref{a_n}) that
\begin{equation}
\label{a_n1}
a_n = \frac{p^{2n-1}-1}{p-1}.
\end{equation}
According to $pa_n+1=(p^{2n}-1)/(p-1)<p^{2n}/(p-1)$, (\ref{M_n}) and (\ref{M_1}), we define $\{M_n\}$ by
\[
M_{n+1} = C_1 p^{-2pn} M_n^{p^2},\quad M_1=C_g \e^p, \quad (C_1:=C_0^{p+1}(p-1)^p)
\]
Hence, we reach to
\[
\begin{array}{ll}
\log M_{n+1}
&=(1+p^2+\cdots+p^{2(n-1)})\log C_1\\
&\quad-\{n+p^2(n-1)+\cdots+p^{2(n-1)}(n-n+1)\}\log (p^{2p})\\
&\quad+p^{2n}\log M_1\\
&\d=\frac{p^{2n}-1}{p^2-1}\log C_1-p^{2n}\sum_{j=1}^n\frac{j}{p^{2j}}\log (p^{2p})+p^{2n}\log M_1\\
&\d\ge-\frac{1}{p^2-1}\log C_1+p^{2n}\left\{\frac{1}{p^2-1}\log C-S_{p^2}\log (p^{2p})+\log M_1\right\},
\end{array}
\]
where
\[
S_{p^2}:=\sum_{j=1}^\infty \frac{j}{p^{2j}} < \infty.
\]
Therefore, we obtain 
\[
U(x,t) \ge \frac{C_1^{-1/(p^2-1)}}{(1+t-x)^{1+a}} \left\{\log\left(\frac{1+t-x}{1+R}\right)\right\}^{-1/(p-1)} \exp \left( K_1(x,t) p^{2n} \right) \ \mbox{in} \ D
\]
where
\[
\begin{array}{ll}
\d K_1(x,t) := &\d \frac{1}{p-1}\log\left\{\log\left(\frac{1+t-x}{1+R}\right)\right\} \\
&\d + \frac{1}{p^2-1}\log C_1-S_{p^2}\log (p^{2p})+\log M_1
\end{array}
\]
If there exists a point $(x_0,t_0)\in D$ such that
\[
K_1(x_0,t_0)>0,
\]
we have
\[
u(x_0,t_0)=\infty
\]
by letting $n\rightarrow\infty$, so that $T<t_0$.
Let us set $2x_0=t_0$ and $4(1+R)^2 < t_0$. Then $K_1(t_0/2,t_0)>0$ follows from
\begin{equation}
\label{con1}
2^{-1}\log\left(t_0\right)C_1^{1/(p+1)}p^{-2p S_{p^2}(p-1)}(C_g \e^p)^{p-1}> 1
\end{equation}
because the inequality 
\[
\log \left(\frac{1+t_0-x_0}{1+R}\right) \ge 2^{-1} \log t_0
\]
hold for $4(1+R)^2 < t_0$.
The condition (\ref{con1}) follows from
\[
\log t_0 > 2\{C_1^{1/(p+1)}p^{-2p S_{p^2}(p-1)}C_g^{p-1}\}^{-1}\e^{-p(p-1)}
\quad\mbox{for}\ 0<\e\le\e_3
\]
if we define $\e_3$ by
\[
4(1+R)^2=2\{C_1^{1/(p+1)}p^{-2p S_{p^2}(p-1)}C_g^{p-1}\}^{-1}\e_3^{-p(p-1)}.
\]The proof for $p(1+a)+b=0$ is now completed.

\par\noindent
{\bf Case 2. $\v{p(1+a)+b<0}$.} 
\par\noindent
In this case, let us retake $(x,t)$ as follows: 
\[
(x,t) \in D_{a,b} := \{ (x,t) \in D \ | \ 1+t-x > 2^{1/(-p(1+a)-b)}(1+R) \} .
\]
By assanption, 
\begin{equation}
\label{halfup}
(1+t-x)^{m(-p(1+a)-b)} - (1+R)^{m(-p(1+a)-b)}> 2^{-1} (1+t-x)^{m(-p(1+a)-b)}
\end{equation}
holds for any $m \ge 1$. As in case 1, assume that an estimate
\begin{equation}
\label{case_p(1+a)+b<0}
U(x,t)\ge \frac{M_n}{(1+t-x)^{1+a}}\left\{(1+t-x)^{-p(1+a)-b}\right\}^{a_n}
\quad\mbox{in}\ D_{a,b}
\end{equation}
and $\{a_n\}$ is defined by (\ref{a_n1}). By appropriately changing the constants of $M_n$ in Case 1, (\ref{case_p(1+a)+b<0}) holds for any $n$.
In fact, by the same calculation as case 1, we have 
\[
\begin{array}{lll}
\d U(x,t) 
& \d \ge C_0\int_{(t+x+R)/2}^t  \frac{|U(x+t-s,s)|^p}{(1+t+x)^{1+a}(1+|2s-t-x|)^{1+b}}ds \\
& \d \ge C_0M_n^p\int_{(t+x+R)/2}^t  \frac{\left\{(1+2s-t-x)^{-p(1+a)-b}\right\}^{pa_n}ds}{(1+t+x)^{1+a}(1+2s-t-x)^{p(1+a)+1+b}}\\
& \d = \frac{C_0M_n^p}{(pa_n+1)(-p(1+a)-b)(1+t+x)^{1+a}} \times \\
& \d \quad \times \left[ \left\{(1+t-x)^{-p(1+a)-b}\right\}^{pa_n+1} - \left\{(1+R)^{-p(1+a)-b}\right\}^{pa_n+1} \right] \\
& \d \ge \frac{C_0M_n^p\left\{(1+t-x)^{-p(1+a)-b}\right\}^{pa_n+1}}{2(pa_n+1)(-p(1+a)-b)(1+t+x)^{1+a}}. \\
\end{array}
\]
and 
\[
\begin{array}{lll}
\d U(x,t) 
& \d\ge  C_0 \int_{(t-x+R)/2}^{t-x}  \frac{|U(x-t+s,s)|^p}{(1+t-x)^{1+a}(1+|2s-t+x|)^{1+b}}ds \\
& \d\ge  \frac{C_0^{p+1}M_n^{p^2}}{2^p(pa_n+1)^p(-p(1+a)-b)^p} \times  \\
& \d \quad \times\int_{(t-x+R)/2}^{t-x}  \frac{\left\{(1+t-x)^{-p(1+a)-b}\right\}^{p^2a_n+p}ds}{(1+t-x)^{1+a}(1+2s-t+x)^{p(1+a)+1+b}}\\
& \d = \frac{C_0^{p+1}M_n^{p^2}\left\{(1+t-x)^{-p(1+a)-b}\right\}^{p^2a_n+p}}{2^p(pa_n+1)^p(-p(1+a)-b)^{p+1}(1+t-x)^{1+a}} \times \\
& \d \quad \times \left[(1+t-x)^{-p(1+a)-b} - (1+R)^{-p(1+a)-b} \right] \\
& \d \ge \frac{C_0^{p+1}M_n^{p^2}\left\{(1+t-x)^{-p(1+a)-b}\right\}^{p^2a_n+p+1}}{2^{p+1}(pa_n+1)^p(-p(1+a)-b)^{p+1}(1+t-x)^{1+a}}.
\end{array}
\]
Therefore, the same argument as case1 holds in case 2 by changing $C_1$ to $C_2 := C_1/2^{p+1}(-p(1+a)-b)^{p+1}$.
we obtain 
\[
U(x,t) \ge \frac{C_2^{-1/(p^2-1)}}{(1+t-x)^{1+a}} \left\{(1+t-x)^{-p(1+a)-b}\right\}^{-1/(p-1)} \exp \left( K_2(x,t) p^{2n} \right) \ \mbox{in} \ D_{a,b}
\]
where
\[
\begin{array}{ll}
\d K_2(x,t) := &\d \frac{1}{p-1}\log\left\{(1+t-x)^{-p(1+a)-b}\right\} \\
&\d + \frac{1}{p^2-1}\log C_2-S_{p^2}\log (p^{2p})+\log M_1
\end{array}
\]
Then $K_2(t_0/2,t_0)>0$ follows from
\begin{equation}
\label{con2}
2^{p(1+a)+b}t_0^{-p(1+a)-b}C_2^{1/(p+1)}p^{-2p S_{p^2}(p-1)}(C_g \e^p)^{p-1}> 1.
\end{equation}
The condition (\ref{con2}) follows from
\[
t_0 > 2\{C_1^{1/(p+1)}p^{-2p S_{p^2}(p-1)}C_g^{p-1}\}^{-1}\e^{-p(p-1)/(-p(1+a)-b)}
\quad\mbox{for}\ 0<\e\le\e_4
\]
if we define $\e_4$ by the condition of $D_{a,b}$ :
\[
2^{1/(-p(1+a)-b))}(1+R)-2^{-1}=\{C_1^{1/(p+1)}p^{-2p S_{p^2}(p-1)}C_g^{p-1}\}^{-1}\e_4^{-p(p-1)/(-p(1+a)-b)}.
\]The proof for $p(1+a)+b<0$ is now completed. $\e_2$ of Proposition \ref{prop:2} should be the smaller one of $\e_3$ and $\e_4$.
\hfill$\Box$

\section*{Acknowledgement}
\par
The author would like to thank Professor Hiroyuki Takamura (Tohoku Univ., Japan)
for his consistent encouragements and discusstions, and Professor Kunio Hidano (Mie Univ. in Japan)
for his pointing out trivial oversights on the regularity of the solution
for low powers of the nonlinear terms.

This work was supported by JST, the establishment of university fellowships towards the creation of science technology innovation, Grant Number JPMJFS2102.


\bibliographystyle{plain}

\end{document}